\documentclass[12pt,twoside,a4paper]{article}
\usepackage{amsmath,amssymb,fancyhdr}
\usepackage{graphicx,color}
\usepackage{float,verbatim}
\usepackage[OT1]{fontenc}
\usepackage{paralist}
\usepackage{tikz}
\usepackage[textwidth=16cm]{geometry}
\usepackage{hyperref}
\usepackage{enumitem}
\newtheorem{thm}{Theorem}[section]
\newtheorem{lem}{Lemma}[section]

\newtheorem{remn}{Remark}[section]
\numberwithin{equation}{section}

\renewcommand{\bar}[1]{\overline{#1}}

\newcommand{\wh}{\widehat}

\newcommand{\wt}{\widetilde}
\renewcommand{\tilde}{\widetilde}

\newcommand{\pf}{{\em Proof. }}
\renewcommand{\Box}{\square}
\newcommand{\qed}{$\quad\,\Box$}
\newcommand{\e}{\mathrm e}

\newcommand{\calf}{\mathcal F}
\newcommand{\calu}{\mathcal U}
\newcommand{\calw}{\mathcal W}

\newcommand{\thmref}[1]{\makebox{Theorem~\ref{#1}}}
\newcommand{\lemref}[1]{\makebox{Lemma~\ref{#1}}}
\newcommand{\secref}[1]{\makebox{Section~\ref{#1}}}
\newcommand{\subsecref}[1]{\makebox{Subsection~\ref{#1}}}

\newcommand{\remref}[1]{\makebox{Remark~\ref{#1}}}

\newcommand{\rr}{{\mathbb R}}

\newcommand{\zz}{{\mathbb Z}}

\newcommand{\const}{\text{const}}

\newcommand{\essinf}{\operatornamewithlimits{ess\,inf}}
\newcommand{\esssup}{\operatornamewithlimits{ess\,sup}}
\newcommand{\norm}[2]{\|{#1}\|_{#2}}

\author{Annegret Glitzky\footnote{Weierstrass Institute for Applied Analysis and Stochastics, Mohrenstr. 39, 10117 Berlin, Germany,
\{glitzky,liero\}@wias-berlin.de} 
\and Matthias Liero}

\title{Uniqueness and regularity of weak solutions of a drift-diffusion system for perovskite solar cells}      

\newcommand{\SUBMISSION}{}

\begin{document}

\maketitle

\begin{abstract}
\noindent
We establish a novel uniqueness result for an instationary drift-diffusion 
model for perovskite solar cells. This model for vacancy-assisted charge transport 
uses Fermi--Dirac statistics for electrons and holes and Blakemore statistics for 
the mobile ionic vacancies in the perovskite. 
Existence of weak solutions and their boundedness was proven in a previous work.   
For the uniqueness proof, we establish improved integrability of the gradients 
of the charge-carrier densities. 
Based on estimates obtained in the previous paper, 
we consider 
suitably regularized
continuity equations 
with partly frozen arguments 
and apply 
the regularity results for scalar quasilinear elliptic equations by  Meinlschmidt \& Rehberg, 
Evolution Equations and Control Theory, 2016, 5(1):147-184. 
\end{abstract}
\ifdefined\SUBMISSION

\fi
\section{Introduction}\label{s:1}
Perovskite materials have outstanding optical and electronic properties \cite{Schmidt2021RoadMap}.
They have the advantages of adjustable band gap, high absorption coefficient, long exciton diffusion length, excellent carrier mobility and low exciton binding energy. Using layers made from the family of halide perovskites in solar cell concepts has shown a considerable potential for high performance and low production costs in solar cells. In recent years, the power conversion efficiency of perovskite solar cells in the laboratory has raised rapidly from 3.8\,\% in 2009 to 25.5\,\% in 2021, see \cite{ZHANG2022421}. In 2023, the highest efficiency perovskite-silicon tandem solar cell had a power conversion efficiency of 29.1\%, see also \cite{NREL2024}. 
However, a number of challenges like stability and durability issues remain before they can become a competitive commercial technology.

The crystal structure of perovskite materials is of the form
ABX$_3$, where A and B represent cations and X is an anion.
Typically, these ions can move and  vacancies remain in the lattice, but their mobility is much lower compared to the mobility of electrons and holes. 
As one aspect, the diffusion engineering of ions and charge carriers for stable efficient perovskite solar cells turns out to be an important task \cite{ag-bi17,kress2022persistent}.

In \cite{ag-vag21}, a drift-diffusion model for the vacancy-assisted charge transport in perovskite solar cells was derived. It describes the heterostructure composed by classical semiconductor materials and perovskite materials. It consists of continuity equations for electrons, holes, and different kinds of vacancies in the perovskite material that are self-consistently coupled to a Poisson equation. Using Fermi--Dirac statistics for electrons and holes and Blakemore statistics for the vacancies, to prevent an unrealistic accumulation of the latter, is motivated in \cite{ag-vag21}. We refer to \secref{sec:model} for a brief summary of the model equations. 
Numerical analysis for a finite-volume discretized version of the model is presented in \cite{ag-abd22}. Moreover, analytical investigations concerning the existence and boundedness of weak solutions of the instationary problem are performed in \cite{ag-gli24}. 
A related model for memristor devices, restricted to the setting of Boltzmann statistics and for one fixed domain, is treated in  \cite{ag-jou23}.

The aim of the present paper is a uniqueness and regularity result for weak solutions of this drift-diffusion system for perovskite solar cells under additional assumptions (see (A5) below)
based on the results in \cite{ag-gli24}. The outline of the paper is as follows.
\secref{sec:model} introduces the drift-diffusion model for 
the vacancy-assisted charge transport in perovskite solar cells.
In \secref{sec:ana}, we formulate the assumptions for its analytical treatment, give the weak formulation of the model equations,  and briefly summarize already established analytical results concerning energy estimates, a priori estimates and existence results. Our main results, i.e., the uniqueness result \thmref{thm:einzig} 
and the higher regularity of the solution (\thmref{higherreg}),
are demonstrated in \secref{sec:unique}. The proof is based on an auxiliary regularity result for scalar quasilinear parabolic equations that we present in \subsecref{meinlschmidt}. In \subsecref{application}, we apply this result to verify the higher regularity result \thmref{higherreg} for solutions of our coupled drift-diffusion system that finally enables us to prove the uniqueness result in \subsecref{einzig}.
Concluding remarks are collected in \secref{sec:conclusion}.
\section{Drift-diffusion modeling of perovskite solar cells}\label{sec:model}

\subsection{Drift-diffusion system}

We consider a rescaled version of a drift-diffusion model derived in \cite{ag-vag21} describing the vacancy-assisted charge transport in perovskite solar cells.
Let $\Omega \subset \rr^d$ be the spatial domain of the solar cell and $I$ denote the index set of mobile carriers.
Additionally to the movement of electrons and holes in $\Omega $, we take into account the migration of ionic vacancies  from the index set $I_0 \subset I$ in $\Omega_0 \subset \Omega$.
We denote the densities of electrons, holes, and ionic vacancies by  $u_i$, $i\in I := \{ n, p \} \cup I_0$, where $i = n$ and $i = p$ correspond to electrons and holes, respectively.
The considered drift-diffusion model couples a Poisson equation for the electrostatic potential $\psi$ self-consistently to the continuity equations for the densities $u_i$
\begin{subequations}\label{eq:scaled}
  \begin{align} \label{eq:Poisson}
    -\nabla \cdot (\varepsilon\nabla \psi) & =
    \begin{cases}
    C+z_nu_n+z_pu_p&\text{  in }(0,\infty)\times (\Omega\setminus \overline{\Omega_0}),  \\
    C+z_nu_n+z_pu_p+\sum_{i\in I_0}z_iu_i&\text{  in }(0,\infty)\times\Omega_0,
    \end{cases}\\
    \frac{\partial u_i}{\partial t}- \nabla\cdot (z_i\mu_iu_i\nabla\varphi_i)
    & = G-R,\,\, i=n,p,\hspace{7.51em} \text{in }(0,\infty)\times\Omega, \label{eq:ElectronHoles} \\
    \label{eq:IonicSpecies}
    \frac{\partial u_i}{\partial t}- \nabla\cdot (z_i\mu_iu_i\nabla\varphi_i)
    &= 0,\quad i\in I_0,\hspace{9.972em}\text{in }(0,\infty)\times\Omega_0,
  \end{align}
\end{subequations}
  where $\varepsilon$ denotes a rescaled dielectric permittivity, $z_i$ the charge number of a species $i \in I$ and $C$ the fixed doping density.
Additionally, $\mu_i$ are the rescaled carrier mobilities and
the generation/recombination terms $G$ and $R$ in the continuity equations for electrons and holes \eqref{eq:ElectronHoles} are discussed in \subsecref{subsec:reac}.
For the ionic vacancies $i \in I_0$  we do not consider any reactions.
The statistical relation connecting the potentials $\varphi_i$ and $\psi$ to
the charge-carrier densities $u_i$ is given by
\begin{align} \label{eq:state-eq}
  u_i =N_i\calf_i(z_i(\varphi_i-\psi)+\zeta_i) =N_i\calf_i(v_i+\zeta_i),
  \quad\text{where}\quad v_i =z_i(\varphi_i-\psi),\quad i\in I,
\end{align}
with the effective densities of state $N_n$ and $N_p$ for electrons and holes, the maximal density of vacancies $N_i$, $i\in I_0$, the chemical potentials $v_i$ and $\zeta_i:=z_i\mathrm{E}_i$, $i\in I$, where $\mathrm{E}_i$ is the band-edge energy.
The statistics function $\mathcal{F}_i$ will be discussed in \subsecref{sec:stat-function}.
Note that in comparison to the model in \cite{ag-vag21}, we rescaled the electrostatic potential $\psi$ and the quasi Fermi potentials $\varphi_i$ by the thermal voltage $U_T = k_B T /q$
($k_B$ is Boltzmann's constant, $T$ the (constant) temperature and $q$  the elementary charge).
Moreover, the chemical potentials $v_i$ and the band-edge energies $\mathrm{E}_i$ are rescaled by $k_B T$.
And we multiplied the dielectric permittivity $\varepsilon$ by $U_T / q$ and the mobilities $\mu_i$ by $U_T$. This rescaling was already used in \cite{ag-gli24}, where the existence of weak solutions of the system \eqref{eq:scaled} was proved.

\subsection{Statistical functions} \label{sec:stat-function}
The statistical functions for electrons and holes in classical (inorganic) semiconductors are given by the Fermi--Dirac integral of order $1/2$ (see e.g., \cite{ag-sze07})
\begin{align} \label{eq:FD-1-2}
  F_{1/2}(z)=\frac{2}{\sqrt{\pi}}\int_0^\infty\frac{\xi^{1/2}}{\exp(\xi-z)+1}\,\text d\xi, \quad\text{for }z\in\rr
\end{align}
i.e., $\mathcal{F}_n = \mathcal{F}_p = F_{1/2}$.
For small to moderate carrier densities, the Fermi--Dirac integral of order $1/2$ can be approximated by the exponential function (Boltzmann statistics), i.e., $F_{1/2}(z)\approx \e^z$, see \cite{ag-sze07}.
Our analysis and the results in \cite{ag-gli24} are not restricted to these specific choices, but work under the following general assumptions on the statistical functions
\begin{equation}\label{eq:F12}
  \left\{
  \begin{aligned}
    (\text{i})\;\, & \mathcal{F}_i\in C^1(\rr), \quad
    \lim_{z\to -\infty} \mathcal{F}_i(z)=0, \quad \lim_{z\to +\infty} \mathcal{F}_i (z)=+\infty,\\
    (\text{ii})\;\, & z\le c(1+\mathcal{F}_i(z)) \quad \text{for}\; z\in\rr_+,\\
    (\text{iii}) \;\,& 0<\mathcal{F}_i'(z)\le \mathcal{F}_i(z)\le \e^z\quad\text{for }z\in\rr,
\end{aligned}
\right.
\quad i = n, p.
\end{equation}
The accumulation of too many vacancies is physically unrealistic, as it would destroy the crystal structure.
Limiting adequately the vacancy concentration can be done via a proper choice of the statistical function.
In \cite{ag-vag21}, the use of the Fermi--Dirac integral of order $-1$ (which corresponds to Blakemore statistics $F_{B,\gamma}$ function with $\gamma = 1$)
\begin{equation}\label{eq:Blakemore}
F_{-1}(z)= F_{B,1}(z),\quad\text{where}\quad F_{B,\gamma}(z) = \frac{1}{\e^{-z}+\gamma} \quad\text{for }z\in\rr,
\end{equation}
i.e.,  $\calf_i = F_{-1}$ for all $i \in I_0$ is motivated and proposed.
In our analytical investigations, we assume that the statistics function of the ionic vacancies satisfies the following properties
\begin{equation} \label{eigf-1}
  \left\{
  \begin{aligned}
    (\text{i})\;\,  & \mathcal{F}_i\in C^2(\rr), \quad\lim_{z\to -\infty} \mathcal{F}_i(z)=0, \quad \lim_{z\to +\infty} \mathcal{F}_i (z)=1,\\
    (\text{ii})\;\,  &\calf_i'(z) <\calf_i(z) <\e^z\quad\text{for }z\in\rr.,\\
    (\text{iii})\;\, & \calf''_{i}(z)<0,\quad
    \frac{|\calf_{i}''(z)|}{\calf'_{i}(z)}<1\quad \text{for}\; z\in\rr_+,\\
    (\text{iv})\;\,  & 1< (\e^z \calf'_{i}(z))^{-1} < c, \quad \text{for}\; z\in\rr_+,
\end{aligned}
\right.
\quad i \in I_0.
\end{equation}
Note that the Fermi--Dirac integral of order $1/2$ and the Boltzmann statistics satisfy the properties \eqref{eq:F12} while the Fermi--Dirac integral of order $-1$ satisfies \eqref{eigf-1}, see also
 \cite{ag-gli24}.

\subsection{Generation-recombination rate and photogeneration}\label{subsec:reac}
According to \cite{ag-far17}, we use for the generation-recombination term $R$
in \eqref{eq:ElectronHoles} the expression 
\begin{equation}\label{eq:reak}
R=r(\cdot,u_n,u_p)\big(1-\e^{\varphi_n-\varphi_p}\big)
\quad\text{with}\quad r(\cdot,u_n,u_p)=r_0(\cdot,u_n,u_p)\,u_nu_p.
\end{equation}
Following \cite{ag-vag21}, the function $r$ is given by the sum of all recombination processes relevant in photovoltaics like radiative and trap-assisted  Shockley-Read-Hall recombination.
In the simplest case, the photogeneration rate $G$ is considered 
to be constant in time and one assumes a Beer-Lambert generation 
profile in the vertical direction $x_\mathrm{vert}$, i.e., $G(x)=F_\mathrm{ph}\alpha_G \,\e^{-\alpha_G x_\mathrm{vert}}$
for $x=(\bar x,x_\mathrm{vert})$, 
where $F_\mathrm{ph}$ denotes the incident photon flux and $\alpha_G$ a material absorption coefficient.

\subsection{Initial and boundary conditions}
We prescribe for all densities  initial values 
\begin{equation}\label{eq:inin}
u_i(0)=u_i^0 \quad \text{ in }\Omega,\quad i=n,p,\quad u_i(0)=u_i^0
\quad\text{ in }\Omega_0,\quad i\in I_0.
\end{equation}
For the formulation of boundary conditions
we decompose $\partial \Omega$ into the set of Ohmic contacts 
$\Gamma_D$ 
and the semiconductor-insulator interface $\Gamma_{\!\!N}$.
Ohmic contacts like semiconductor-metal interfaces are modeled by Dirichlet boundary
conditions 
\begin{subequations}\label{RB}
  \begin{equation}\label{RB1}
  \psi=\psi^D,\quad
  \varphi_i= \varphi^D,\quad i=n,p,\quad
  \text{on }(0,\infty)\times \Gamma_D.
  \end{equation}
Semiconductor-insulator interfaces are described by no-flux boundary conditions
\begin{equation}\label{RB2}
\varepsilon \nabla \psi\cdot \nu = \mu_iu_i\nabla \varphi_i\cdot \nu  
 = 0,\quad i=n,p,\quad\text{on }(0,\infty)\times \Gamma_N,
\end{equation}
where $\nu$ denotes the outer normal vector. 
At the boundary of the perovskite domain $\partial\Omega_0$ with  outer normal vector $\nu_0$ we suppose no normal flux of ionic vacancies 
\begin{equation}
 \mu_iu_i\nabla \varphi_i\cdot\nu_0  = 0\quad\text{on }(0,\infty)\times \partial\Omega_0,\,\,i\in I_0.
\end{equation}
\end{subequations}

\section{Analysis of the instationary drift-diffusion model}
\label{sec:ana}

\subsection{Notation and assumptions on the data}

 In our estimates,  
positive constants, depending at most on the data of our problem,
are denoted by $c$. In particular, we allow them to change from line to line.  
We work with the Lebesgue spaces $L^p(\Omega)$ and the
Sobolev spaces $W^{1,p}(\Omega)$, $p\in[1,\infty]$, 
and $H^1(\Omega)=W^{1,2}(\Omega)$. 
For $p\in[1,\infty]$, we define the spaces $W^{1,p}_D(\Omega)$
as the closure of the set 
$
\{y|_\Omega: y\in C_0^\infty(\rr^d), \text{ supp\,} y\cap \Gamma_D
=\emptyset\}
$
in the Sobolev space $W^{1,p}(\Omega)$ and $W^{-1,p}_D(\Omega):=W^{1,p'}_D(\Omega)^*$, where
$1/p+1/p'=1$.

We study the 
drift-diffusion model under the assumptions
\begin{enumerate}[label=(A\arabic*),series=l_after]
\item
$\Omega, \Omega_0\subset\rr^2$ are bounded Lipschitz domains, 
$\Omega\cup\Gamma_N$ is regular in the sense of \\
Gr\"oger \cite{ag-gro89},
$\Gamma_D,\,\Gamma_N\subset\Gamma:=\partial \Omega$ are disjoint subsets
 such that
$\overline{\Gamma_D\cup\Gamma_N}=\Gamma$\\ 
and $\Gamma_D$ is closed and
$\text{mes}(\Gamma_D)>0$,
$\Omega_0\subset\Omega$, 
$\Omega_n=\Omega_p:=\Omega$, $\Omega_i:=\Omega_0$, $i\in I_0$. 
\item
$\calf_i$ fulfill \eqref{eq:F12} for $i=n,p$, 
$\calf_i$ fulfill \eqref{eigf-1} for $i\in I_0$,\\
$N_i,\,\mu_i\in L^\infty(\Omega_i)$,  
$0<\underline{N}\le N_i\le \overline{N},\,
0<\underline{\mu}\le \mu_{i}\le \overline{\mu}$
a.e. in $\Omega_i$, $\zeta_i=\const$,  
 $i\in I$,\\
and $C,\,\varepsilon\in L^\infty(\Omega)$, 
$0<\underline{\varepsilon}\le\varepsilon$ a.e. in $\Omega$,
$z_n=-1$, $z_p=1$, $z_i\in\zz$, $i\in I_0$,\\
$v_0^D:=\psi^D,\,\varphi^D\in 
W^{1,\infty}(\Omega)$.  
\item 
$G\in L^\infty(\rr_+;L^\infty_+(\Omega))$,\\
$R=r(\cdot,u_n,u_p)\big(1-\e^{\varphi_n-\varphi_p}\big)$,
such that $r(\cdot,u_n,u_p)= r_0(\cdot,u_n,u_p)\,u_n\,u_p$, where\\
$r_0:\Omega\times[0,+\infty)^2
\to \rr$ is a
 Carath\'eodory function with\\
$
0\le r_0(\cdot,u_n,u_p)\le \overline{r}$
for all 
 $(u_n,u_p)\in [0,+\infty)^2\text{ and a.a.\ }x\in \Omega$.
\item  
$u_i^0\in L^\infty(\Omega)$, $0<\underline{u}\le u_i^0\le \bar u $ 
a.e.\ in $\Omega$, $i=n,p$,
\\
 $u_i^0\in L^\infty(\Omega_0)$, $0<\underline{u}\le u_i^0\le \bar{u_i} < N_i$ a.e.\ in $\Omega_0$, $i\in I_0$.
\end{enumerate}

In the following, we suppress in the writing the spatial 
position $x$ in the terms $r$ and $r_0$.
 

\subsection{Weak formulation}\label{subsec:weak}

We introduce the functions $e_i:\rr\to(0,\infty)$, $i=n,p$, and
$e_i:\rr\to(0,1)$, $i\in I_0$, for the statistical relations 
\begin{equation}\label{eq:estat}
e_i(y)=\calf_i(y+\zeta_i),\quad i\in I.
\end{equation}
The inverse functions $e_i^{-1}$ are well-defined on 
$(0,\infty)$ for $i=n,p$, and  on $(0,1)$ for $i\in I_0$.
From Assumption (A2), we obtain the following relation between the densities $u_i$
and the associated chemical potentials $v_i$
 \begin{equation}\label{eq:gradient}
 \begin{split}
&u_i=N_ie_i(v_i)=N_i\calf_i(v_i+\zeta_i),\quad
v_i=e_i^{-1}(\tfrac{u_i}{N_i})
=\calf_i^{-1}(\tfrac{u_i}{N_i})-\zeta_i,
\\
&
\nabla\tfrac{u_i}{N_i}=\calf_i'(v_i+\zeta_i)\nabla v_i=
e_i'(v_i)\nabla v_i.
\end{split}
\end{equation}
We define the function spaces
\begin{align*}
V&:=V_D^3\times V_0^{\# I_0},\quad
V_D:=\{y\in H^1(\Omega): y|_{\Gamma_D}=0\},\quad
V_0=H^1(\Omega_0),\\[0.5em]
H&:=V_D\times L^2(\Omega)^2\times L^2(\Omega_0)^{\# I_0},\quad
Z:=H^1(\Omega)\times L^\infty(\Omega)^2\times L^\infty(\Omega_0)^{\# I_0},
\\[0.5em]
U&:=\big\{u\in V_D^*\times L^\infty(\Omega)^2\times L^2(\Omega_0)^{\# I_0}:
\ln u_i\in L^\infty(\Omega),\,i=n,p,\\
&\hspace{2cm}\,0<\essinf_{x\in\Omega_0} u_i/N_i \le
\esssup_{x\in\Omega_0} u_i/N_i <1,\,\,i\in I_0\big\},
\end{align*}
where $\# I_0$ denotes the number of elements of the set $I_0$.

In \cite{ag-gli24}, we considered a weak formulation of \eqref{eq:scaled} in 
the form 
$u'+A(v,v)=0$, $u=E(v)$, $u(0)=u^0
$
 with the variables 
$v=(v_0,v_n,v_p,(v_i)_{i\in I_0}):=(\psi,
(z_i(\varphi_i{-}\psi))_{i\in I} )$ (potentials) and
$u:=(u_0,u_n,u_p,(u_i)_{i\in I_0})$ (densities)
as well as operators $A$ and $E$ defined below (see also \cite{ag-gaj89,ag-gaj96a,ag-gli18}).
The first component $u_0$ represents the total charge density of the device, i.e.,
the right-hand side in \eqref{eq:Poisson}, and the initial value 
$u^0:=(u_0^0,u_n^0,u_p^0,(u_i^0)_{i\in I_0})$ is such that
$\langle u_0^0,w\rangle_{V_D}=\sum_{i\in I}\int_{\Omega_i} z_i u_i^0w\,\mathrm{d}x+\int_\Omega Cw\,\text dx$ for all 
$w\in V_D$. 
We recall that $z_n=-1$, $z_p=1$, and $z_i$, $i \in I_0$, stands for the charge number of the 
 $i$-th vacancy species. 
In these variables, the problem reads
\begin{equation}\label{eq:vscaled}
\begin{split}
&-\nabla \cdot (\varepsilon\nabla v_0)  = 
\begin{cases} 
C+z_nu_n+z_pu_p, &\text{  in }(0,\infty)\times(\Omega\setminus\overline{\Omega_0})\\
C+z_nu_n+z_pu_p+\sum_{i\in I_0}z_iu_i, &\text{  in }(0,\infty)\times\Omega_0
\end{cases},\\
& \frac{\partial u_i}{\partial t}- \nabla\cdot 
(\mu_iu_i(\nabla v_i+z_i\nabla v_0))  = G-R,
\quad
u_i=N_ie_i
(v_i),\quad i=n,p,\\
& \frac{\partial u_i}{\partial t}
- \nabla\cdot (\mu_iu_i(\nabla v_i+z_i\nabla v_0))  = 0,\quad 
u_i=N_ie_i
(v_i), \quad i\in I_0,  
\end{split}
\end{equation}
where 
$R=r(u_n,u_p)\big(1-\e^{-v_n-v_p}\big)$.
The prescribed Dirichlet values are
$v_n^D:=v_0^D-\varphi^D$, $v_p^D:=\varphi^D-v_0^D=-v_n^D$.
In case of the ions, we do not have to prescribe a Dirichlet value, however,
for a unified notation, we set $v^D:=(v_0^D,v_n^D,v_p^D,(0)_{i\in I_0})$.
We introduce the operators 
$E_0:v_0^D+V_D\to V_D^*$,
$E:(v^D+V)\cap Z\to V^*$, and $A:Z\times(v^D+V)\to V^*$ by
\begin{equation*}
\begin{split}
E(v) & :=(E_0(v_0),(N_ie_i(v_i))_{i\in I}),
\quad \langle E_0 (v_0),\bar v_0\rangle_{V_D}  :=
\int_\Omega\varepsilon\nabla v_0\cdot\nabla\bar v_0\,\mathrm{d}x,
\\
\langle A (w,v),\bar v \rangle_V  &:=
\sum_{i\in I}\int_{\Omega_i} N_ie_i(w_i)\mu_i\nabla(v_i+z_iv_0)\cdot\nabla(\bar v_i+z_i\bar v_0)
\,
\mathrm{d}x\\
&\quad\quad+\int_\Omega  \,[r(N_ne_n(w_n),N_pe_p(w_p))
\big(1-\e^{-w_n-w_p}\big)-G]
(\bar v_n+\bar v_p)\,\mathrm{d}x
\end{split}
\end{equation*}
for all $\bar v\in V$.
The weak formulation of the drift-diffusion system \eqref{eq:scaled}, 
\eqref{eq:inin}, \eqref{RB} 
is
\begin{equation}\label{eq:weak}
\tag{P}
  \begin{split}
& u'+A(v,v)=0,\quad u=E(v)\quad \text{a.e.\  on }\rr_+,\quad u(0)=u^0,\\
& u\in H^1_{\text{loc}}(\rr_+,V^*),\quad v-v^D\in L^2_{\text{loc}}(\rr_+,V)
\cap L^\infty_{\text{loc}}(\rr_+,Z).
\end{split}
\end{equation}%

\subsection{Summary of analytical results for the drift-diffusion system}
Next we give a very short overview on analytical results for 
Problem \eqref{eq:weak} obtained so far in \cite{ag-gli24}. 
They concern energy estimates, solvability, and bounds of solutions.

{\bf Energy estimates.}
In the analysis of the drift-diffusion problem, entropy methods as in
\cite{ag-gaj89,ag-gaj96a,ag-jue16} play an important part. In \cite{ag-gli24}, we considered a free energy functional $\Psi$
containing an electrostatic contribution and a chemical part that is  obtained from the statistical relations of the different types of species. For states $u=Ev\in V^*\cap U$, it has the form
\begin{equation}\label{eq:energy}
\begin{split}
\Psi(u) 
&=\int_\Omega\frac{\varepsilon}{2}|\nabla (v_0{-}v_0^D)|^2\,\mathrm{d}x
+\sum_{i\in I} 
\int_{\Omega_i} \int_{v_i^D}^{v_i}[u_i{-}N_ie_i(y)]\mathrm{d}y\,
\mathrm{d}x.
\end{split}
\end{equation}
For Blakemore statistics, 
$e_i=F_{-1}(\cdot+\zeta_i)=F_{B,1}(\cdot+\zeta_i)$, the chemical energy for an ionic vacancy $i\in I_0$ in 
\eqref{eq:energy} reads 
\[
\int_{\Omega_0} \int_0^{v_i}[u_i-N_ie_i(y)]\mathrm{d}y\,
\mathrm{d}x
=\int_{\Omega_0}
\Big(u_i\ln \frac{u_i}{N_i}+(N_i-u_i)\ln\big(1-\frac{u_i}{N_i}\big)
-u_i\zeta_i+N_i\ln (\e^{\zeta_i}+1)\Big)\,\text dx,
\]
which forces the ion vacancy density $u_i$ to stay in $[0,N_i]$.
For arguments $u=E(v)$ we find
\begin{equation}\label{eq:psiunten}
\norm{u_n}{L^1(\Omega)} +
\norm{u_p}{L^1(\Omega)} +
\sum_{i\in I_0}\norm{u_i}{L^1(\Omega_0)}+\norm{v_0}{H^1(\Omega)}^2
\le c(1+ \Psi(u)).
\end{equation}
In \cite[Theorem~4.1]{ag-gli24} the following energy estimate was obtained. 

\begin{thm}[{\cite[Theorem~4.1]{ag-gli24}}]\label{thm:fe}
Let {\rm (A1) -- (A4)} be fulfilled. Then there is a constant $c>0$ such that
\[
\Psi(u(t))\le \big(\Psi(u(0))+c\big)\e^{ct}\quad \forall \,\,t>0
\] 
for any solution $(u,v)$ to  Problem 
\eqref{eq:weak}. Additionally, if the Dirichlet values are compatible with 
thermodynamic equilibrium (meaning $v_0^D,\,\varphi^D=\const$) and if the generation rate $G$ is identically zero, then  
$t\mapsto\Psi(u(t))$ is monotonically decreasing.
\end{thm}

Moreover, for any solution $(u,v)$ to \eqref{eq:weak} we have the conservation laws 
 (see \cite[Remark~4.1]{ag-gli24})
\[
\int_{\Omega_0}u_i(t)\,\mathrm{d} x
= \int_{\Omega_0}u_i^0\,\mathrm{d} x,\quad 
i\in I_0,\quad\text{for all }t\in\rr_+,\quad \forall i\in I_0.
\]

{\bf Existence and boundedness of solutions.}
Under the Assumptions {\rm (A1) -- (A4)}
the existence of a solution to 
Problem~\eqref{eq:weak} is verified by demonstrating the existence of solutions for any finite time interval
 $S:=[0,T]$. We introduce the Problem
\begin{equation}\label{ps}
\tag{P$_{\text S}$}
  \begin{split}
&u'+A(v,v)  =  0, \quad u=E(v)\,\text{ a.e.\ on }S,\quad     
u(0)  =  u^0,\\
&u\in H^1(S,V^*),\quad 
v-v^D \in L^2(S,V)\cap L^\infty(S,Z).
\end{split}
\end{equation}%

\begin{thm}[{\cite[Theorem~5.1]{ag-gli24}}]\label{thm:exist}
We assume {\rm (A1) -- (A4)}.
Then, for all $T>0$, $S:=[0,T]$, there exists at least one solution to Problem
\eqref{ps}.
\end{thm} 
The proof of this existence result is based on the following steps: 
First, a regularized problem (P$_{\text M}$) 
on the time interval $S$ is discussed, where the
state equations as well as the reaction term are regularized 
(with parameter $M$). Solvability of (P$_{\text M}$) is shown
by time discretization, derivation of suitable a priori estimates, and passage to the limit 
(see \cite[Lemma~5.1 and Lemma~5.2]{ag-gli24}). 

Then, we provide
a priori estimates for solutions $(u^M,v^M)$ to (P$_{\text M}$) that are 
independent of $M$ (see \cite[Lemma~5.6]{ag-gli24}, here we use Moser 
techniques to get positive lower bounds for the carrier densities $u_i^M$.
Moreover, in \cite[Lemma~5.5 and Lemma~5.7]{ag-gli24}  upper bounds independent of $M$  are derived 
for the chemical potentials $v_i^M$, $i\in I$).
Thus, a solution to (P$_{\text M}$) is a solution to (P$_{\text S}$), 
if $M$ is chosen sufficiently large.

\begin{thm}[{\cite[Theorem~5.2]{ag-gli24}}]\label{thm:globounds}
We assume {\rm (A1) -- (A4)}.
Then, for all $T>0$, $S=[0,T]$ there exist $c_0(T),\,c_1(T)$, $c_2(T)>0$ 
such that for any solution $(u,v)$ to Problem \eqref{ps} 
\begin{equation*}
\begin{split}
& c_1(T)\le u_i(t)\le c_0(T)\quad  \text{a.e. in }\Omega,\quad i=n,p,\\ 
& c_1(T)\le u_i(t)\le c_2(T)N_i\quad  \text{a.e. in }\Omega_0,\quad i\in I_0,\quad \forall \,t\in S.
\end{split}
\end{equation*}
\end{thm}

\begin{remn}
In \cite[Theorem~5.2]{ag-gli24}, the upper bounds for electron and hole densities were obtained by test functions of the form 
\[
2^m\big(0,[(\tfrac{u_n}{N_n}-K)^+]^{2^m-1},
[(\tfrac{u_p}{N_p}-K)^+]^{2^m-1},(0)_{i\in I_0}\big), \quad m\ge 1, 
\]
where
$
K:=\max\big\{\max_{i=n,p}
\norm{e_{i}(z_i(\varphi^D-v_0^D))}{L^\infty(\Omega)},\max_{i\in I}
\norm{u_i^0/N_i}{L^\infty(\Omega_i)}\big\}
$
and Moser-type estimates.  
For $i\in I$, the positive lower bound for $u_i$ was established by test functions
\[
-2^m\frac{\big[\big(\ln{(u_i/N_i)}+\wh K\big)^-\big]^{2^m{-}1}}{u_i/N_i}
,\quad m\ge 1,\text{ where}
\]
$
\wh K:=\max\big\{\max_{i=n,p}
\norm{\ln e_{i}(v_i^D)}{L^\infty},\max_{i\in I}
\norm{\big(\ln (u_i^0/N_i)\big)^-}{L^\infty},
\max_{i\in I}\ln e_{i}(0)\big\}
$ used separately for the $i$-th continuity equation
and Moser techniques.
Finally, the upper bounds for the densities of ionic vacancies strictly below $N_i$ were (roughly speaking) obtained by applying the last line of \eqref{eigf-1} and
by means of test functions
\[
\frac{2^m\,\big[(\e^{v_i}-\wt K)^+\big]^{2^m{-}1} \,\e^{v_i}}{e_i'(v_i)}
,\quad m\ge 1,\quad\text{where}\quad
\wt K:=\max_{i\in I_0}\max\big\{\e^{\norm{e_i^{-1}(u_i^0/N_i)}{L^\infty(\Omega_0)}},\e^{-\zeta_i},1\big\}
\]
used separately for the $i$-th continuity equation together with a Moser iteration.
\end{remn}


\section{Uniqueness result and 
regularity properties}\label{sec:unique}

One main result of this paper is the 
uniqueness of the weak solution to the drift-diffusion system for perovskite solar cells in Problem \eqref{eq:weak}.
Note that for the pure electronic charge transport (only van Roosbroeck system for electrons and holes and $I_0=\emptyset$) the following results are known:
In the case of Boltzmann statistics, uniqueness of solutions holds (see e.g.\ \cite[Theorem 5.1]{ag-gaj90}).
Note that in the case of Boltzmann statistics and for homogeneous materials with constant coefficients, 
the diffusive term in the charge-carrier flux has the form $\mu_i\nabla u_i$ whereas
for more general statistical relations terms of the form $b_i(u_i)\nabla u_i$ have to be handled to show uniqueness.
For Fermi--Dirac statistics uniqueness has been proven under unjustified ad hoc assumptions concerning the smoothness of solutions in \cite[Theorem 3.2]{ag-gaj89}. 
Another uniqueness result for a special coupled system of one nonlinear parabolic equation and one elliptic equation with more general statistical relations  is derived in \cite{ag-gaj03}.

In addition to the Assumptions (A1) -- (A4), needed already for the 
existence proof and for the derivation of the boundedness 
results of \thmref{thm:globounds}, we will suppose the following properties of the data
for our uniqueness result:
\begin{enumerate}[label=(A\arabic*),resume*=l_after]
  \item
$\Omega$ and $\Omega_0$ are domains with Lipschitz boundary\\ (strong Lipschitz domains
in the sense of Grisvard
\cite[Chapter 1.2]{ag-gri85});\\
$N_i=\const$, $\calf_i\in C^2(\rr)$, $i\in I$;\\
$u_i^D:=N_ie_i(v_i^D)$,  
$u_i^0-u_i^D\in W^{1,\lambda}_D(\Omega)$, $i=n,p$,\\ 
$u_i^0\in W^{1,\lambda}_D(\Omega_0)$, $i\in I_0$, for some $\lambda>2$;\\
$r_0$ is locally Lipschitz continuous.
\end{enumerate}
We recall that the subset of $\partial\Omega_i$, where Dirichlet boundary conditions are prescribed in the continuity equations, is empty for the ionic vacancies $i\in I_0$, 
we  only set for a unified writing $W^{1,\lambda}_D(\Omega_i):=W^{1,\lambda}(\Omega_i)$ and correspondingly 
$u_i^D:=0$ for $i\in I_0$.

Now we are able to formulate the main result of our present paper that concerns the unique solvability of Problem (P).
 
\begin{thm}\label{thm:einzig}
Under the  Assumptions {\rm (A1) -- (A5)} the Problem {\rm (P)} has at most one solution.
\end{thm}

\pf 
For the proof of \thmref{thm:einzig}, it suffices to demonstrate that for any $T>0$ and $S:=[0,T]$, we have the uniqueness of the solution to \eqref{ps} (see \thmref{thm:einzigs}). Namely, if we would have two different solutions $(u^1,v^1)$ and $(u^2,v^2)$ to Problem (P) then there would be some $T_0>0$ and $S_0=[0,T_0]$ such that $(u^1|_{[0,T_0]},v^1|_{[0,T_0]})$ and $(u^2|_{[0,T_0]},v^2|_{[0,T_0]})$ would be two different solutions to (P$_{\rm S_0}$). Therefore the desired result follows directly from \thmref{thm:einzigs} below.
\qed

\begin{thm}\label{thm:einzigs}
We assume {\rm (A1) -- (A5)}. Let $T>0$ be arbitrarily given and $S:=[0,T]$. Then Problem \eqref{ps} has at most one solution.
\end{thm}

We start with a short outline of the proof of \thmref{thm:einzigs}:

Step 1. Since $N_i$ is constant (see (A5)), we obtain from \eqref{eq:gradient} for solutions $(u,v)$ to \eqref{ps} that $\nabla u_i=N_ie_i'(v_i)\nabla v_i$ and the diffusive part in the flux terms of the continuity equations in \eqref{eq:vscaled} can be rewritten as
\[
\mu_iu_i\nabla v_i
= \mu_i\frac{u_i}{N_i}(e_i^{-1})'\big(\frac{u_i}{N_i}\big)\nabla u_i= b_i(u_i)\nabla u_i,\quad\text{where }\, b_i(u_i):=\mu_i\frac{u_i}{N_i}(e_i^{-1})'\big(\frac{u_i}{N_i}\big),\quad i\in I.
\]

According to the uniform positive lower and upper bounds for the concentrations $u_i$ guaranteed by \thmref{thm:globounds}, we find constants $\underline {b}$ and  $\overline{b}$ such that
\[
0<\underline {b}\le b_i(u_i(t))\le\overline{b} \quad \text{a.e. in }
\Omega_i\quad\forall t\in S,\quad i\in I,
\]
for all solutions $(u,v)$ to \eqref{ps}. 
Under Assumptions (A2) and (A5), 
for the given $v_i^D\in W^{1,\infty}(\Omega)$ and constant $N_i$, we define $u_i^D:=N_ie_i(v_i^D)\in W^{1,\infty}(\Omega)$ and obtain
$u_i=\wt u_i+u_i^D\in L^2(S,H^1(\Omega))\cap L^\infty(S,L^\infty(\Omega))$ with 
$\wt u_i:=N_ie_i(v_i)-N_ie_i(v_i^D)\in L^2(S,W^{1,2}_D(\Omega))$,
$i=n,p$. (Note that for the ionic vacancies we do not prescribe any Dirichlet conditions. Only for a unified notation we set $u_i^D=0$, $i\in I_0$.)
Thus, solutions $(u,v)$ to \eqref{ps} fulfill
\begin{align*}
\int_0^T\langle u_i',w_i\rangle_{V_D}\,\text dt
&=-
\int_0^T\!\!\int_\Omega \big((b_i(u_i)\nabla u_i{+}z_i\mu_iu_i\nabla v_0)\cdot\nabla w_i +Q(u_n,u_p)w_i\big)\,\text dx\,\text dt,
\quad \! i=n,p,\\
\int_0^T\langle u_i',w_i\rangle_{V_0}\,\text dt&=-
\int_0^T\!\!\int_{\Omega_0} (b_i(u_i)\nabla u_i+z_i\mu_iu_i\nabla v_0)\cdot\nabla w_i\, \text dx\,\text dt,
\quad i\in I_0,
\end{align*}
for all $w_i\in L^2(S,W^{1,2}_D(\Omega))$, $i=n,p$, and all
$w_i\in L^2(S,H^1(\Omega_0))$, $i\in I_0$, respectively, with
\begin{equation}\label{eq:defQ}
Q(u_n,u_p):=r(u_n,u_p)\big(1-\e^{-e_n^{-1}(u_n/N_n)-\e_p^{-1}(u_p/N_p)}\big)-G.
\end{equation}
In this setting, we will apply for each separate continuity equation the regularity theory for quasilinear parabolic equations, see e.g., \cite{ag-mei16,ag-mei17,ag-mei23}.

Step 2. 
The aim is to improve the integrability properties of $\nabla u_i$.
We verify the existence of a $q>2$ such that for all solutions $(u,v)$ to \eqref{ps} and all $s\ge 1$ we find $u_i\in L^s(S,W^{1,q}(\Omega_i))$, $i\in I$, see \thmref{higherreg} in \subsecref{application}.

For this purpose, we exploit the regularity result \cite[Theorem 5.3]{ag-mei16} 
for scalar quasilinear parabolic equations, which is summarized in \subsecref{meinlschmidt}. 
We apply this result in \subsecref{application}, separately to
suitably regularized versions of the above continuity equations, $i\in I$, 
where $u_i$ in $b_i$ and in the drift term are replaced by truncated densities. The cut-off is chosen in such a way that 
the truncations have no effect for solutions due to Theorem~\ref{thm:globounds}.
Moreover, in these regularized equations, the electrostatic potential $v_0$ in 
the drift term and the carrier densities $u_n,\,u_p$ 
in the generation-recombination term $Q(u_n,u_p)$ are fixed by a solution $(u,v)$ to \eqref{ps}. 
\thmref{thm:parreg} below guarantees a unique solution
$y_i \in L^s(S,W^{1,q}_D(\Omega))$ ($i=n,p$) or 
$y_i \in L^s(S,W^{1,q}(\Omega_0))$ ($i\in I_0$) to the regularized continuity equation (see Problem \eqref{piq} below). 
This solution is also a solution in $H^1(S,W^{1,2}_D(\Omega_i)^*)\cap L^2(S,W^{1,2}_D(\Omega_i))$ 
to an auxiliary Problem \eqref{pi2}, see 
\subsecref{application}. 

However, for any solution $(u,v)$  to \eqref{ps} 
(where $u_i^D=0$ for $i\in I_0$) the function $\wt u_i = u_i-u_i^D$ 
is also a solution to the auxiliary Problem \eqref{pi2}, and we establish that $y_i=\wt u_i$
has to hold (see \lemref{lem:same}). Thus, we obtain the better integrability property of $\nabla u_i$,
 which proves the assertion of \thmref{higherreg} in \subsecref{application}.

Step 3. Exploiting this improved regularity of $\nabla u_i$, 
we establish in \subsecref{einzig} the uniqueness result stated in \thmref{thm:einzigs} for the coupled system.

\subsection{Higher regularity result for scalar quasilinear parabolic
 equations}\label{meinlschmidt}

Let us first recall a well-known regularity result for the linear elliptic situation
due to Gröger. Let $\rho:\Omega\to \rr$
be a bounded, measurable function with 
$0<\underline\rho\le\rho(x)\le \bar \rho $ a.e. in $\Omega$. 
We define the linear operator, denoted in the following by  
$-\nabla\cdot\rho\nabla +1: \,W^{1,2}_D(\Omega)\to W^{-1,2}_D(\Omega)$, via
\[
\langle (-\nabla\cdot\rho\nabla+1)w,\bar w\rangle_{W^{1,2}_D(\Omega)}
:=
\int_\Omega(\rho\nabla w\cdot\nabla \bar w +w\bar w)\,\text dx
\quad \text{ for }w,\bar w\in W^{1,2}_D(\Omega).
\]
We also denote the maximal restriction of $-\nabla\cdot\rho\nabla+1$ to any 
of the spaces $W^{-1,q}_D(\Omega)$ ($q>2$) by the same symbol,
$-\nabla\cdot\rho\nabla +1$. 

If the set $\Omega\cup\Gamma_N$ is regular in the sense of Gr\"oger, 
the regularity result \cite[Theorem 1]{ag-gro89} ensures an exponent $q_0>2$ 
such that 
$-\nabla\cdot\rho\nabla +1$ maps $W^{1,q}_D(\Omega)$ onto 
$W^{-1,q}_D(\Omega)$ for all $q\in[2,q_0]$. In particular, 
the operator is a topological isomorphism, and its inverse 
is Lipschitz continuous. We emphasize that the exponent $q_0$ 
only depends on the geometrical setting $\Omega\subset\rr^d,\,\Gamma_N$ 
as well as $\bar\rho$, and $\underline\rho$.

Next, we summarize the regularity result for scalar quasilinear parabolic equations, 
obtained in \cite{ag-mei16}, that we will apply to improve the integrability 
properties of the gradient of the charge-carrier densities $u_i$ for solutions to \eqref{ps}. 
We consider equations of the following form
\begin{equation}
  \label{eq:ScaQuasPara}
  y'(t)-\nabla\cdot \big(\theta(y(t))\mu\nabla y(t)\big) +y(t)=
\calf(t,y(t)),\quad y(0)=y^0.
\end{equation}

The regularity result covers the geometric setting characterized by the following assumption.

{\bf Assumption (A)} 
(see Assumptions 2.2, 2.4 in \cite{ag-mei16})\\
Let $\Omega\subset\rr^d$ be a bounded domain and $\Gamma_D$ be a 
closed subset of $\partial\Omega$ (which may be empty), and
$\Gamma_N:=\partial\Omega\setminus \Gamma_D$ such that\\
i) If $x\in \partial\Omega\setminus\bar{\Gamma_N}$, there 
is a domain $\calu_x=:\calu$ with $x\in\calu$, such that 
$\calu\cap\Gamma_N=\emptyset$ and $\calu\cap\Omega$ has 
only finitely many connected components $Y_1,\dots,Y_k$, 
where $x$ is a limit point of each $Y_j$.
Moreover, for every $j\in\{1,\dots,k\}$, there exists a
$\tau_j>0$, an open neighbourhood $U_j$ of $x$ 
satisfying $Y_j\subseteq U_j\subseteq\calu$, and a 
bi-Lipschitz mapping $\phi_j$ defined on an open 
neighbourhood of $\bar{U_j}$ into $\rr^d$,
such that $\phi_j(x)=0$,
$\phi_j(U_j)=\tau_j\{x\in \rr^d:\norm{x}{\infty}<1\}$,
$\phi_j(Y_j)=\tau_j\,\{x\in \rr^d:\norm{x}{\infty}<1,\,x_d<0\}$, 
and 
$\phi_j(\partial Y_j\cap U_j)=\tau_j\,\{x\in \rr^d:
\norm{x}{\infty}<1, \,x_d=0\}$.\\[1ex]
ii) For each $x\in \bar{\Gamma_N}$ there is an open 
neighbourhood $U_x=:U$ of $x$, a number $\tau_x=:\tau>0$ and 
a bi-Lipschitz mapping $\phi_x=:\phi$ from an open  
neighbourhood of $\bar U$ into $\rr^d$, such that 
$\phi(x)=0\in\rr^d$, 
$\phi(U)=\tau \{x\in \rr^d:\norm{x}{\infty}<1\}$,
$\phi(U\cap\Omega)=\tau \{x\in \rr^d:\norm{x}{\infty}<1,
\,x_d<0\}$, and 
$\phi(\partial\Omega\cap U)=\tau \{x\in \rr^d:\norm{x}{\infty}
<1,\,x_d=0\}$.\\
\hspace*{1cm} (a)
If $x\in \Gamma_N$, then $U$ does not intersect $\Gamma_D$, i.e.,
$\phi(\Gamma_D\cap U)=\emptyset$.
\\
\hspace*{1cm} (b)
If $x\in \bar{\Gamma_N}\cap\Gamma_D$, then $\phi(\Gamma_D\cap U)=
\tau\{x\in \rr^d:\norm{x}{\infty}<1,\,x_d=0,\,x_{d-1}\le 0\}$.
\\[1ex]
iii) Each of the occurring mappings $\phi$ is volume-preserving.

\begin{remn}
In \cite[Definition 1.3.12]{ag-mei17} a similar condition on the geometric setting as in Assumption (A) is introduced, called volume-preserving generalized regular in the sense of Gr\"oger. (But note that, Gr\"oger regularity in the original work \cite{ag-gro89} is formulated in terms of $\Omega\cup\Gamma_N$).
\end{remn}

Furthermore, the assumptions concerning the right-hand side $\calf$ 
in \eqref{eq:ScaQuasPara} are as follows.

{\bf Assumption (B)} 
(see Assumption 5.1 in \cite{ag-mei16})\\
 i) The function $\calf:S\times C(\bar\Omega)
 \to W^{-1,q}_D(\Omega)$ is a Carath\'eodory function, i.e., 
 $\calf(\cdot,y):S\to
 W^{-1,q}_D(\Omega)$ is measurable for each $y\in C(\bar\Omega)$ and 
$\calf(t,\cdot):C(\bar\Omega)\to W^{-1,q}_D(\Omega)$
is continuous for all $t\in S$.\\[1ex]
ii) For $s\in (1,\infty)$, we assume that the superposition operator $y\mapsto [t\mapsto \calf(t,y(t))]$ is continuous from every bounded set of 
$C(\overline{S\times \Omega})$ to $L^s(S,W^{-1,q}_D(\Omega))$
with 
\[
\sup_{y\in C(\bar\Omega)}
\norm{\calf(\cdot,y)}{L^s(S,W^{-1,q}_D(\Omega))}\le c_\calf
\]
for some constant $c_\calf>0$. 

\begin{remn}\label{alt}
Note that the continuity condition in Assumption (B) is satisfied for a Carath\'eodory function $\calf$, if the boundedness property is fulfilled and for every $\eta>0$ there exists a function $L_\eta\in L^s(S)$ such that
\[
\norm{\calf(t,y_1)-\calf(t,y_2)}{W^{-1,q}_D(\Omega)}\le 
L_\eta(t)\norm{y_1-y_2}{C(\bar\Omega)}\quad\text{f.a.a. }t\in S,
\]
where $y_1,\,y_2\in C(\bar\Omega)$ with
$\norm{y_1}{C(\bar\Omega)},\,\norm{y_2}{C(\bar\Omega)}\le \eta$.
\end{remn}

\begin{thm}(Theorem 5.3 in \cite{ag-mei16})\label{thm:parreg}
We assume {\rm (A)} with $d\le 3$.
Let $\mu$ be a measurable coefficient function on $\Omega$ with
$\underline\mu\le \mu\le \overline \mu$ a.e. in $\Omega$. 
We suppose that $\theta:\rr\to [\underline{\theta},\overline{\theta}]$  with 
$0<\underline{\theta}<\overline{\theta}$ is Lipschitz continuous on bounded sets. Assume further that for some $q>d$ the map
$-\nabla\cdot\mu\nabla+1:W^{1,q}_D(\Omega)\to
W^{-1,q}_D(\Omega)$ is a topological isomorphism and let 
$s>2(1-\tfrac{d}{q})^{-1}$ be such that the initial value $y^0$ 
is an element in the real interpolation space $(W^{1,q}_D(\Omega),W^{-1,q}_D(\Omega))_{1/s,s}$.
Let $\calf:S\times C(\bar\Omega)\to W^{-1,q}_D(\Omega)$ satisfy the Assumption {\rm (B)} for this $s$. 
Then there exists a global solution $y\in W^{1,s}(S, W^{-1,q}_D(\Omega))\cap
L^s(S,W^{1,q}_D(\Omega))$ of the quasilinear equation  \eqref{eq:ScaQuasPara}.
If $\calf$ additionally fulfils the assumptions in \remref{alt} then
this solution is unique.
\end{thm}

Let us mention that Theorem 5.3 in \cite{ag-mei16} is originally formulated 
for the space dimension $d=3$. A consultation with the authors of the corresponding paper, 
J.~Rehberg and H.~Meinlschmidt,  confirmed the validity of their result also for space dimensions $d\le 3$.

Indeed, the necessity of $d\le 3$ is a bit hidden, it is needed to guarantee 
uniformity of the domains of each of the operators
$y\mapsto -\nabla\cdot\theta(y)\mu\nabla  +1$ for $y\in W_D^{1,q}(\Omega)$, 
which is only available for space dimensions up to 3, see also 
\cite[Lemma 5.5]{ag-mei16} and \cite[Lemma 6.2]{ag-dis15}.

For completeness, we briefly discuss this fact in the two-dimensional setting.
The essential point is ensuring for localized problems of the equation 
$-\nabla\cdot(\vartheta\mu\nabla y)+y=f$ with $f\in W^{-1,q}_D(\Omega)$
for some $q\in[2,\infty)$ that the new right-hand sides $f_j$ also belong to 
$W^{-1,q}_D(\Omega)$ with the same exponent $q$, provided 
that $\vartheta:\Omega\to\rr$ is a uniformly continuous 
function with 
$0<\underline\theta\le \vartheta\le \overline \theta$ 
on $\Omega$.
More precisely, let us consider 
a finite covering $\calw_1,\,\dots,\calw_l$ of $\Omega$ and a 
subordinate partition of unity $\eta_1,\dots,\eta_l$ on $\Omega$,
such that the localized equations read
\begin{equation}\label{eq:sub}
-\nabla\cdot \big(\vartheta\mu\nabla(\eta_j y)\big)+\eta_j y=f_j,\quad j=1,\dots,l,
\end{equation}
where $f_j$ denotes the new right-hand side due to localization via $\eta_j$.
If $f_j\in W^{-1,q}_{D}(\Omega)$, we can follow the proof of \cite[Lemma 6.2]{ag-dis15}. 
Then the small variation of the coefficient function 
$\vartheta$ on each of the $\calw_j$, perturbation 
arguments, and Gr\"oger's localization technique ensure 
the uniformity of the domains of each of the operators
$-\nabla\cdot\theta(y)\mu\nabla y +y$ needed in 
\cite[Lemma 5.5]{ag-mei16}.

We argue as follows:
For $y\in W^{1,q}_D(\Omega)$ and 
$\eta_j\in W^{1,\infty}(\Omega)$ with support in $\calw_j$, we find 
$\eta_j y\in W^{1,q}_D(\Omega)$.
For $C_0^\infty$ test functions $\phi$  and 
 $\eta_j\in W^{1,\infty}_0(\calw_j)$ we derive
\begin{equation*}
\begin{split}
\int_\Omega \vartheta\mu\nabla(\eta_j y)\cdot\nabla \phi\,\text dx
 = & - \int_\Omega \phi\vartheta\mu\nabla y \cdot \nabla \eta_j\,\text dx
+\int_\Omega y\vartheta\mu\nabla \eta_j \cdot\nabla\phi\,\text dx\\
& +\int_\Omega \vartheta\mu\nabla y \cdot \nabla(\eta_j\phi)\,\text dx.
\end{split}
\end{equation*}
Let $\Omega_\bullet:=\Omega\cap \calw_j$ and
$\Gamma_{D_\bullet}:=\Gamma_D\cap\calw_j$. We use the notation 
$W^{1,q}_{D_\bullet}(\Omega_\bullet)$, which means 
$W^{1,q}(\Omega_\bullet)$ if $\Gamma_D\cap\calw_j=\emptyset$. 
Then $f_j$ in \eqref{eq:sub} has the form
\[
f_j=-\mu\vartheta\nabla y|_{\Omega_\bullet}\cdot\nabla \eta_j|_{\Omega_\bullet}
+T_y+f_{\eta_i}\in W^{-1,q}_{D_\bullet}(\Omega_\bullet),
\]
where $f_{\eta_i}: w\mapsto \langle f, \wt{\eta_jw}\rangle _{H^1_D(\Omega)}$ (the tilde denotes the prolongation of $\eta_jw$ by zero on $\Omega$) and $T_y$ denotes the form
$w\mapsto\int_{\Omega_\bullet}y\mu\vartheta \nabla \eta_j\cdot \nabla w\,\text dx$.
Note that if $\Omega_\bullet$ is a Lipschitz domain we apply the Sobolev
embedding  theorems, thus the term $\mu\vartheta\nabla y\cdot\nabla \eta_j$ is generically in $L^2$ and can be interpreted as an element of $W^{-1,q}_{D_\bullet}(\Omega_\bullet)$ for $q\in[2,\infty)$. In summary we obtain  $f_j\in W^{-1,q}_{D_\bullet}(\Omega_\bullet)$, and $f_j\in W^{-1,q}_{D}(\Omega)$.

\begin{remn}
A statement for quasilinear equations clearer formulated for space dimension $d=2$ and $d=3$ can be found in the dissertation of H. Meinlschmidt \cite[Theorem 2.2.12]{ag-mei17}. For the domain and the position of the boundary conditions it is assumed that $\Omega\cup\Gamma_D$ is volume-preserving generalized regular in the sense of Gr\"oger, compare  \cite[Definition 1.3.12]{ag-mei17}. 
\end{remn}

\begin{remn}
Recently, the result of \cite[Theorem 5.3]{ag-mei16} was extended 
in \cite[Theorem 3.1]{ag-mei23} and existence and uniqueness of global-in-time solutions in the $W^{-1,q}_D$ - $W^{1,q}_D$-setting for abstract quasilinear parabolic PDEs with non-smooth data and mixed boundary conditions, including a nonlinear source term
with at most linear growth was demonstrated (for $d=2$ and $d=3$). 
Subsequently, the authors used a bootstrapping argument
to achieve improved regularity, in particular H\"older-continuity, of these global-in-time solutions for 
the abstract equation under suitable additional assumptions. 
\end{remn}
 

\subsection{Application of the regularity result to the regularized  continuity equations}\label{application}

For given $v_i\in L^2(S,H^1(\Omega))\cap L^\infty(S,L^\infty(\Omega))$ and $v_i^D\in W^{1,\infty}(\Omega)$ and constant $N_i$, we define $u_i^D:=N_ie_i(v_i^D)\in W^{1,\infty}(\Omega)$ and obtain
$u_i=\wt u_i+u_i^D\in L^2(S,H^1(\Omega))\cap L^\infty(S,L^\infty(\Omega))$ with 
$\wt u_i:=N_ie_i(v_i)-N_ie_i(v_i^D)\in L^2(S,W^{1,2}_D(\Omega))$,
$i=n,p$. For the ionic vacancies we do not prescribe any Dirichlet conditions. Only for a unified notation we set $u_i^D=0$, $i\in I_0$.

Because of (A1), (A2), \cite[Theorem 1]{ag-gro89} guarantees some $\pi>2$ such that
$-\nabla\cdot\varepsilon\nabla+1:W^{1,q}_D(\Omega)\to
W^{-1,q}_D(\Omega)$ is a topological isomorphism 
for all $q\in [2,\pi]$.

We recall that for solutions to \eqref{ps} it holds that 
$u_i\in C(S,L^2(\Omega_i))$, $i\in I$. 
(This is obtained by $v_i\in L^2(S,W^{1,2}(\Omega_i))\cap L^\infty(S,L^\infty(\Omega_i))$ which leads by (A5) to $\nabla u_i=N_ie_i'(v_i)\nabla v_i\in L^2(S,L^2(\Omega_i))^2$. 
With $u_i^D=N_ie_i(v_i^D)$ we find $u_i-u_i^D\in L^2(S,W^{1,2}_D(\Omega))$, $i=n,p$. 
And from 
$\{u_i-u_i^D\in L^2(S,W^{1,2}_D(\Omega_i)):(u_i-u_i^D)'\in
L^2(S,W^{-1,2}_D(\Omega_i))\}\subset C(S,L^2(\Omega_i))$ for 
$i=n,p$ as well as 
$\{u_i\in L^2(S,W^{1,2}(\Omega_i)):u_i'\in
L^2(S,W^{1,2}(\Omega_i)^*)\}\subset C(S,L^2(\Omega_i))$  
for $i\in I_0$
(see 
\cite[p.\,147]{ag-gaj74})
we find $u_i\in C(S,L^2(\Omega_i))$ for all $i\in I$.)

Using the Lipschitz continuous dependence of $(-\nabla\cdot\varepsilon\nabla +1)^{-1}:W^{-1,\pi}_D(\Omega)\to W^{1,\pi}_D(\Omega)$ on the right-hand side (see \cite[Theorem 1]{ag-gro89}) and the continuous embedding $L^2(\Omega_i)\hookrightarrow W^{-1,\pi}_D(\Omega_i)$ it follows that
$v_0\in C(S,W^{1,\pi}(\Omega))$. 

We rewrite the continuity equations in \eqref{ps} in the setting of \thmref{thm:parreg}.
Let $i\in I$. With the truncated densities (compare \thmref{thm:globounds})
\[
d_i(w):=
\begin{cases}
\max\{c_1(T),\min\{w,c_0(T)\}\},\quad &i=n,p,\\
\max\{c_1(T),\min\{w,c_2(T)N_i\}\},\quad &i\in I_0,
\end{cases}
\]
we introduce the functions 
\begin{equation*}
\begin{split}
\theta_i(y)& :=\frac{d_i(y{+}u_i^D)}{N_i}(e_i^{-1})'\Big(\frac{d_i(y{+}u_i^D)}{N_i}\Big),\quad i=n,p,\\ 
\theta_i(y)& :=\frac{d_i(y)}{N_i}(e_i^{-1})'\Big(\frac{d_i(y)}{N_i}\Big), \quad i\in I_0.
\end{split}
\end{equation*}
Since the quantities $(e_i^{-1})'$ are bounded from below and above for arguments
in $[\tfrac{c_1(T)}{N_i},\tfrac{c_0(T)}{N_i}]$, $i=n,p$,
and in 
$[\tfrac{c_1(T)}{N_i},c_2(T)]$, $i\in I_0$, we find $0<\underline \theta<\overline \theta<\infty$ such that
$\underline \theta\le\theta_i(y)\le\overline \theta$ for all $y\in \rr$, $i\in I$. Moreover, since $e_i\in C^2(\rr)$, the truncated functions 
$\theta_i$ are Lipschitz continuous.

We fix an arbitrary solution $(u,v)$ to \eqref{ps} and use 
the corresponding quantities $v_0(t)$, $u_n(t)$, and $u_p(t)$ to define
right-hand sides $\calf_i:S\times C(\bar\Omega)\to W_D^{-1,q}(\Omega)$ by
\begin{equation*}
\begin{split}
\langle  {\calf}_i(t,y),\bar w\rangle _{W_D^{1,q'}(\Omega)}
:=-\int_\Omega\Big\{&
(\theta_i(y) \mu_i\nabla u_i^D
+\mu_i z_i d_i(y+u_i^D)\nabla v_0(t))\cdot\nabla \bar w\\
&
+Q(u_n(t),u_p(t))\bar w-(u_i(t)-u_i^D)\bar w\Big\}\,\text dx\quad \forall \bar w\in W_D^{1,q'}(\Omega),
\end{split}
\end{equation*}
$i=n,p$, with $Q$ given in \eqref{eq:defQ}. For $i\in I_0$ we introduce $\calf_i: S\times C(\bar\Omega_0)\to W_D^{-1,q}(\Omega_0)$, 
\[
\langle \calf_i(t,y),\bar w\rangle_{W_D^{1,q'}(\Omega_0)} 
:=-\int_{\Omega_0}\Big\{
\mu_iz_i d_i(y)\nabla v_0(t)\cdot\nabla \bar w -u_i(t)\bar w\Big\}
\,\text dx,\quad \forall \bar w\in W_D^{1,q'}(\Omega_0)\quad i\in I_0.
\]

We intend to
apply the result for a scalar quasilinear parabolic equation formulated in 
\thmref{thm:parreg} separately for each of the continuity equations.
We fix the exponent $q>2$ as follows: By (A1), (A2), \cite[Theorem 1]{ag-gro89} ensures some $q^*>2$ such that
$-\nabla\cdot\mu_i\nabla+1:W^{1,q}_D(\Omega_i)\to
W^{-1,q}_D(\Omega_i)$, $i\in I$, are topological isomorphisms 
for all $q\in [2,q^*]$. With $\lambda$ from (A5), we now fix
\begin{equation}\label{eq:exp}
q:=\min\{\lambda,\pi,q^*\}>2.
\end{equation}

Separately, for each $i\in I$, we set in \thmref{thm:parreg} 
\[
\Omega:=\Omega_i, \quad
\calf:=\calf_i, \quad
\theta:=\theta_i, \quad
\mu:=\mu_i, \quad
y^0=u_i^0-u_i^D,\quad
s>\frac{2q}{q-2}.
\]
According to \cite[Remark 2.5 iv)]{ag-mei16} and the additional Assumption (A5),
$\Omega_i$ admits bi-Lipschitzian boundary charts which are volume preserving, 
and the Assumption (A) for \thmref{thm:parreg} is fulfilled in our 
geometrical setting. 
 
Note that for the real interpolation spaces $(W^{1,q}_D(\Omega_i),W^{-1,q}_D(\Omega_i))_{1/s,s}$ we have the relations
\begin{equation*}
\begin{split}
W^{1,q}_D(\Omega_i)=W^{1,q}_D(\Omega_i)\cap W^{-1,q}_D(\Omega_i)
& \subset (W^{1,q}_D(\Omega_i),W^{-1,q}_D(\Omega_i))_{1/s,s}\\
& \subset W^{1,q}_D(\Omega_i) +W^{-1,q}_D(\Omega_i),
\end{split}
\end{equation*}
see \cite[Chapter I,\,2]{ag-ama95} and thus by (A5), the functions $y_i^0=u_i^0-u_i^D$, $i\in I$, are admissible 
initial values for \thmref{thm:parreg}.

Due to the definition of the $\calf_i$ and the estimates for the solution $(u,v)$ to \eqref{ps} from \thmref{thm:globounds} 
(the bounds of the carrier densities and the fact that 
$v_0\in C(S,W^{1,\pi}(\Omega))$) and the regularity of the 
Dirichlet values $u_i^D$, we establish that 
$\calf_i:S\times C(\overline{\Omega_i})\to W^{-1,q}(\Omega_i)$ is a Carath\'eodory function (measurability in $t$ for each $y\in C(\bar\Omega)$ and continuity in $y$, 
\begin{equation}\label{help0}
\norm{\calf_i(t,y_1)-\calf_i(t,y_2)}{W^{-1,q}_D(\Omega_i)}\le 
c(\norm{\nabla u_i^D}{L^q(\Omega_i)}+
\norm{\nabla v_0}{C(S,L^q(\Omega))})\norm{y_1-y_2}{C(\bar\Omega_i)}
\end{equation}
 and for all $y_1,\,y_2\in C(\bar\Omega_i)$ with a constant 
$c>0$ (not depending on $t$) for all $t\in S$),
$i\in I$. Furthermore, for $s\in(1,\infty)$, we can estimate for the superposition operator
\[
\norm{\calf_i(\cdot,y_1(\cdot)){-}\calf_i(\cdot,y_2(\cdot))}{L^s(S,W^{-1,q}_D(\Omega_i))}\le 
c(\norm{\nabla u_i^D}{L^q(\Omega_i)}+
\norm{\nabla v_0}{C(S,L^q(\Omega))})\norm{y_1{-}y_2}{C(\bar{S\times\Omega_i})}
\]
for all $y_1,\,y_2\in C(\bar{S\times\Omega_i})$.

Moreover, by the definition of $\calf_i$ and the a priori estimates for the solution $(u,v)$ to \eqref{ps}, we verify the boundedness result
\begin{equation*}
\begin{split}
&\sup_{y\in C(\bar\Omega)}\norm{\calf_i(\cdot,y)}{L^s(S,W^{-1,q}_D(\Omega_i))}\\
&\le c(s)\big(\bar\theta\bar\mu
\norm{u^D_i}{W^{1,q}(\Omega_i)}+
\bar\mu\max\{c_0(T), C_2(T)\bar N\}\norm{v_0}{L^\infty(S,W^{1,q}(\Omega))}\\
&\hspace{8ex}+
\bar r c_0(T)^2(1+\e^{e_n^{-1}(c_0(T))+e_n^{-1}(c_0(T))})+1
+\norm{u_i-u_i^D}{L^\infty(S,L^\infty(\Omega_i))}\big)
\le \wh c(s,T). 
\end{split}
\end{equation*}
Thus, Assumption (B) is fulfilled for all $s$, especially for $s>\frac{2q}{q-2}$. Therefore, in summary, \thmref{thm:parreg}, \remref{alt} and \eqref{help0} guarantee for all $i\in I$ unique 
solutions $y_i$ of the problems
\begin{equation}\label{piq}
\tag{P$_{iq}$}
  \begin{split}
&y_i\in W^{1,s}(S,W^{-1,q}_D(\Omega_i))
\cap L^s(S,W^{1,q}_D(\Omega_i))\quad\text{such that }\\
& y_i'(t)-\nabla\cdot\big(\theta_i(y_i(t))\mu_i\nabla y_i(t)\big)+y_i(t)
=\calf_i(t,y_i(t)),\quad y_i(0)=y_i^0.
\end{split}
\end{equation}%
Since the time interval $S$ is finite, for our initial values $y_i^0\in W^{1,q}_D(\Omega_i)$ (see (A5) and \eqref{eq:exp}) we obtain the integrability 
$y_i\in W^{1,s}(S,W^{-1,q}_D(\Omega_i))
\cap L^s(S,W^{1,q}_D(\Omega_i))$ for all $s\ge 1$.

On the other hand,  we consider the following auxiliary Problems \eqref{pi2} given in weak form
\begin{equation}\label{pi2}
\tag{P$_{i2}$}
  \begin{split}
& w_i\in H^1(S,W^{1,2}_D(\Omega_i)^*)\cap L^2(S,W^{1,2}_D(\Omega_i))
\text{  such that }\\
&
w_i'(t)+B_i(t,w_i(t))=0\quad\text{a.e. in }S,\quad w_i(0)=y_i^0,
\end{split}
\end{equation}%
where 
\begin{equation*}
\begin{split}
&\langle B_i(t,w_i),\bar w_i\rangle_{V_D}\\
&:=
\int_\Omega  \Big\{\mu_i\theta_i(w_i)\nabla (w_i+u_i^D)\cdot \nabla \bar w_i +(w_i-(u_i(t)-u_i^D))\bar w_i+Q(u_n(t),u_p(t))\bar w_i\\
&\qquad\qquad +\mu_iz_i d_i(w_i+u_i^D)
\nabla v_0(t)\cdot\nabla \bar w_i
\Big\}\,\text dx,\quad i=n,p,\\
&\langle B_i(t,w_i),\bar w_i\rangle_{V_0}\\
&:=
\int_{\Omega_0} \Big\{\mu_i\theta_i(w_i)\nabla w_i\cdot \!\nabla \bar w_i +(w_i-u_i(t))\bar w_i
+\mu_iz_i d_i(w_i)
\nabla v_0(t)\cdot\!\nabla \bar w_i
\Big\}\,\text dx,
\end{split}
\end{equation*}
$i\in I_0$, with $Q$ from \eqref{eq:defQ}
for the fixed solution $(u,v)$ to \eqref{ps}.

Then the solution $y_i$ to \eqref{piq} is also a solution 
to Problem \eqref{pi2}, $i\in I$. Moreover, due to Assumption (A5) and 
\thmref{thm:globounds} the function 
$\wt u_i=u_i-u_i^D$ is a solution to Problem \eqref{pi2}, too.

\begin{lem}\label{lem:same}
   We assume {\rm (A1) -- (A5)}. Let $y_i$ denote the unique solution of $y_i'-\nabla\cdot(\theta_i(y_i)\mu_i\nabla y_i)
  +y_i =\calf_i(y_i)$ and $\tilde u_i = u_i-u^D$ with $u_i$ solving the continuity equation
  in \eqref{ps}. Then, $y_i=\tilde u_i$ holds.
\end{lem}

\pf
We test the equations in \eqref{pi2} for the solutions $\wt u_i$ and $y_i$ by $\wt u_i-y_i\in L^2(S, W^{1,2}_D(\Omega_i))$, respectively, 
subtract them, and obtain for $i\in I$
\begin{equation*}
\begin{split}
&\frac{1}{2}\norm{\wt u_i(t){-}y_i(t)}{L^2(\Omega_i)}^2\\
&=
\int_0^t\int_{\Omega_i}\Big\{
{-}\mu_i\theta(\wt u_i)|\nabla(\wt u_i{-}y_i)|^2-(\wt u_i{-}y_i)^2\\
& \hspace{2.0cm}-\mu_i(\theta_i(\wt u_i){-}\theta_i(y_i))\nabla (y_i+u_i^D)\cdot
\nabla (\wt u_i{-}y_i)\\
& \hspace{2.1cm}{-}\mu_iz_i(d_i(\wt u_i{+}u_i^D)-d_i(y_i+u_i^D))\nabla v_0\cdot
\nabla(\wt u_i-y_i)\Big\}\,\text dx\,\text ds\\
&\le \int_0^t\Big\{
-\min\{1,\underline{\mu}\underline{\theta}\}
\norm{\wt u_i-y_i}{H^1(\Omega_i)}^2\\
&\hspace{1.47cm}+c\norm{\wt u_i{-}y_i}{L^\beta(\Omega_i)}\norm{\nabla(y_i{+}u_i^D)}{L^q(\Omega_i)}
\norm{\wt u_i{-}y_i}{H^1(\Omega_i)}\\
& \hspace{1.47cm}+c\norm{\wt u_i{-}y_i}{L^\beta(\Omega_i)}\norm{\nabla v_0}{L^q(\Omega)}
\norm{\wt u_i{-}y_i}{H^1(\Omega_i)}\Big\}\,\text ds\\
&\le \int_0^t\Big\{
-\min\{1,\underline{\mu}\underline{\theta}\}
\norm{\wt u_i{-}y_i}{H^1(\Omega_i)}^2\\
& \hspace{1.47cm}+c\norm{\wt u_i{-}y_i}{L^2(\Omega_i)}^{2/\beta}
\norm{\wt u_i{-}y_i}{H^1(\Omega_i)}^{2-2/\beta}
\big(\norm{\nabla(y_i{+}u_i^D)}{L^q(\Omega_i)}
+\norm{\nabla v_0}{L^q(\Omega)}\big)\Big\}\,\text ds\\
&\le c\int_0^t\Big\{\norm{\wt u_i{-}y_i}{L^2(\Omega_i)}^2\big(
\norm{\nabla(y_i{+}u_i^D)}{L^q(\Omega_i)}^\beta
+\norm{\nabla v_0}{L^q(\Omega)}^\beta\big)\Big\}\,\text ds \quad\forall t\in S,
\end{split}
\end{equation*}
where $1/\beta+1/q=1/2$. In the last but one estimate we used Gagliardo--Nirenberg's inequality 
$\norm{w}{L^\beta(\Omega_i)}\le c \norm{w}{L^2(\Omega_i)}^{2/\beta}
\norm{w}{H^1(\Omega_i)}^{1-2/\beta}$,
and in the last one we applied Young's inequality.
According to \thmref{thm:parreg} and (A5), $y_i+u_i^D\in L^\beta(S,W^{1,q}(\Omega_i))$, $i\in I$. Moreover, \thmref{thm:globounds} and the definition of $q$ in \eqref{eq:exp} ensure that $v_0\in L^\beta(S,W^{1,q}(\Omega))$. Therefore, we obtain by Gronwall's lemma that $\wt u_i=y_i$, $i\in I$. 
\qed

Concluding, \lemref{lem:same} yields the higher regularity of $u_i=\wt u_i+u_i^D\in 
L^s(S,W^{1,q}(\Omega_i))$ for all $s\ge 1$, $i\in I$, for the original solution $(u,v)$ to Problem \eqref{ps}. 
Note that due to our choice the exponent $q$ is uniform for all possible solutions $(u,v)$ to Problem \eqref{ps}.
The arguments in \subsecref{application} have proven the following higher integrability result.

\begin{thm}\label{higherreg}
We assume {\rm (A1) -- (A5)}. Then there exists an exponent $q>2$  such that
\[
v_0\in L^s(S,W^{1,q}_D(\Omega)),\quad 
u_i\in L^s(S,W^{1,q}_D(\Omega_i))\cap W^{1,s}(S,W^{-1,q}_D(\Omega_i)),
\quad i\in I,
\]
for all $s\ge 1$  
and for any solution $(u,v)$ to Problem \eqref{ps}. 
(Note  our definition $W^{1,q}_D(\Omega_0)=W^{1,q}(\Omega_0)$ for all $i\in I_0$.)
\end{thm}

\subsection{Uniqueness of solutions to {\bf (P$_\textbf{S}$)}, 
proof of \thmref{thm:einzigs}}\label{einzig}

\emph{Proof of \thmref{thm:einzigs}}\\
\emph{Step 1.} As already mentioned at the beginning of \secref{sec:unique}, under Assumption (A5) also the charge-carrier densities $ u_i=N_ie_i(v_i)$ are functions in $L^2(S,H^1(\Omega_i))$, $i\in I$. 
Moreover, due to \thmref{thm:globounds},
for solutions $(u,v)$ to \eqref{ps}
the functions  $b_i(u_i):=\mu_i\tfrac{u_i}{N_i}(e_i^{-1})'(\tfrac{u_i}{N_i})$, $i\in I$,
are well-defined and fulfil 
$0<\underline b\le b_i(u_i(t))\le \overline b$ a.e. in $\Omega_i$ for all $t\in S$, $i\in I$. 
Since $e_i\in C^2(\rr)$ and $e_i(z)>0$, we obtain
the local Lipschitz continuity of $(e_i^{-1})'$. This property then
gives the Lipschitz continuity  for the bounded intervals of densities. Similarly, by (A5) we get the Lipschitz continuity of $Q$ defined in \eqref{eq:defQ} for the bounded intervals of densities.

\emph{Step 2.} Let $(u,v)$ and $(\wh u,\wh v)$ be two solutions to Problem \eqref{ps}. We define the differences $\bar u:=u-\wh u$ and $\bar v:=v-\wh v$. Then,
from the Poisson equations we obtain
$
E_0(v_0(t))-E_0(\wh v_0(t))=\bar u_0(t)$ f.a.a. $t\in S$
and 
\begin{equation}\label{diffpg}
\norm{\bar v_0(t)}{V_D}\le c\sum_{i\in I}\norm{\bar u_i(t)}
{L^2(\Omega_i)},
\end{equation}
where $c$ is independent of $t\in S$.

\emph{Step 3.}
Testing $u'+A(v,v)=0$ and $\wh u\,'+A(\wh v,\wh v)=0$ by 
$(0,(\bar u_i)_{i\in I})\in L^2(S,V)$ and using $Q$ from \eqref{eq:defQ},
we find  the estimate
\begin{equation*}
\begin{split}
\frac{1}{2}\sum_{i\in I}\norm{\bar u_i(t)}{L^2(\Omega_i)}^2
&=\int_0^t\Big\{\sum_{i\in I}\int_{\Omega_i}\big\{{-}b_i(u_i)|\nabla\bar u_i|^2\\
& \hspace{2cm} -
\big(\big(b_i(u_i){-}b_i(\wh u_i)\big)\nabla \wh u_i
+\mu_i\bar u_iz_i\nabla v_0+\mu_i \wh u_i z_i\nabla \bar v_0\big)
\cdot\nabla\bar u_i\big\}\text dx\\
& \quad
+\int_\Omega\big(Q(\wh u)-Q(u)\big
)(\bar u_n+\bar u_p)\, \text dx\Big\}\,
\text ds\\
&\le 
\sum_{i\in I}\int_0^t\int_{\Omega_i}
\Big\{{-}\underline b|\nabla \bar u_i|^2+c|\bar u_i||\nabla\bar u_i|
(|\nabla \wh u_i|+|\nabla v_0|)\\
&\hspace{3cm}
+c \wh u_i|\nabla\bar v_0||\nabla\bar u_i| +c|\bar u_i|^2\Big\}\,\text dx\,\text ds\qquad \forall t\in S. 
\end{split}
\end{equation*}
Exploiting $q>2$ from \thmref{higherreg} and setting $s=\beta$ such that $1/s=1/2-1/q$, \eqref{diffpg}, and the Gagliardo--Nirenberg inequality
$\norm{y}{L^s}\le c\norm{y}{L^2}^{2/s}\norm{y}{H^1}^{1-2/s}$ 
for $y\in H^1(\Omega_i)$,
 we continue the estimate via
\begin{equation*}
\begin{split}
\frac{1}{2}\sum_{i\in I}\norm{\bar u_i(t)}{L^2(\Omega_i)}^2
& \le 
\sum_{i\in I}\int_0^t\Big\{
-\underline b \norm{\bar u_i}{H^1}^2
+c\norm{\bar u_i}{L^s} \norm{\bar u_i}{H^1}
(\norm{\nabla \wh u_i}{L^q}+\norm{\nabla v_0}{L^q})\\
&\hspace{2.5cm}
+c\norm{\nabla \bar v_0}{L^2} \norm{\bar u_i}{H^1}
+c\norm{\bar u_i}{L^2}^2\Big\}\,\text ds\\
& \le 
\sum_{i\in I}\int_0^t\Big\{
-\underline b \norm{\bar u_i}{H^1}^2
+c\norm{\bar u_i}{L^s} \norm{\bar u_i}{H^1}
(\norm{\nabla \wh u_i}{L^q}+1)\\
&\hspace{2.5cm}
+c\sum_{j\in I}\norm{\bar u_j}{L^2}\norm{\bar u_i}{H^1}
+c\norm{\bar u_i}{L^2}^2\Big\}\,\text ds\\
& \le 
\sum_{i\in I}\int_0^t\Big\{
-\underline b \norm{\bar u_i}{H^1}^2
+c\norm{\bar u_i}{L^2}^{2/s}\norm{\bar u_i}{H^1}^{2-2/s}
(\norm{\nabla \wh u_i}{L^q}+1)\\
&\hspace{2.5cm}
+c\sum_{j\in I}\norm{\bar u_j}{L^2}\norm{\bar u_i}{H^1}
+c\norm{\bar u_i}{L^2}^2\Big\}\,\text ds.
\end{split}
\end{equation*}
We apply Young's inequality and establish the inequality
\[
\sum_{i\in I}\norm{\bar u_i(t)}{L^2(\Omega_i)}^2
\le 
\sum_{i\in I}\int_0^t c \norm{\bar u_i(t)}{L^2(\Omega_i)}^2
(\norm{\nabla \wh u_i}{L^q(\Omega_i)}^s + 1)\,\text ds.
\]
Since $\wh u_i\in L^s(S,W^{1,q}(\Omega_i))$ for $i\in I$ by \thmref{higherreg}, 
we obtain by Gronwall's lemma that $\bar u_i=0$ meaning 
$u_i=\wh u_i$, $i\in I$. Finally, the inequality \eqref{diffpg} guarantees for the electrostatic potentials $v_0=\wh v_0$. Therefore, also $u_0=\wh u_0$, $v_i=\wh v_i$ for $i\in I$  and 
$(u,v)=(\wh u,\wh v)$ holds. This finishes the proof of \thmref{thm:einzigs}.
\qed

\section{Concluding remarks}
\label{sec:conclusion}
In this paper, we studied a vacancy-assisted charge transport model from an analytical point of view. We demonstrated a new uniqueness result for a drift-diffusion model for perovskite solar cells. We proved the uniqueness of weak solutions under the (restricting) Assumptions (A1) -- (A5). For this purpose, under the additional Assumption (A5) a novel higher integrability and regularity result for weak solutions of the perovskite solar cell model was established. Let us collect some remarks:

1. A corresponding result can also be obtained for the more 
generalized situation, where the  different ionic vacancies live on different subdomains 
$\Omega_i\subset\Omega$ of the domain $\Omega$
with $\Omega_i\neq \Omega_j$ for $i\neq j\in I_0$. The necessary results concerning existence and boundedness of solutions can also be obtained for this situation, see \cite[Section~6]{ag-gli24}.
For the uniqueness result, we have to suppose that $\Omega, \,\Omega_i$, $i\in I_0$, are domains with Lipschitz boundary in Assumption (A5).
The result remains especially true if $\Omega_i=\Omega$ for all $i\in I$.
Then the setting of the paper on memristor devices \cite{ag-jou23} in two space dimensions is covered by our result. However, there the situation in three space dimensions is handled with only one type of ionic species (which is crucial for the analysis therein) with Boltzmann statistics for electrons, holes and the ionic species but without generation/recombination and photogeneration. Whereas for this system no uniqueness result is given in this paper, for
the fast-relaxation limit (solutions of the stationary continuity equations for electrons and holes are substituted in the right-hand side of the Poisson equation) in two space dimensions existence and weak-strong uniqueness analysis is performed in \cite{ag-jou23}.

2.
In \cite[p.~248]{ag-gaj03} the condition that the function $x\mapsto e_i''(x)/e_i'(x)$ is non-increasing for the statistical relation $e_i$ plays an important role in the derivation of a uniqueness result for a special coupled system consisting of a nonlinear parabolic equation and an elliptic equation. Note that the Fermi-Dirac statistics (see e.g.\ the introduction in \cite{ag-gaj03}) as well as Blakemore statistics $F_{B,1}=F_{-1}$ possess this property, in particular
\[
\frac{F_{-1}''(x)}{F_{-1}'(x)}
=\frac{\e^{-x}(\e^{-x}-1)}{(\e^{-x}+1)^3}
\frac{(\e^{-x}+1)^2}{\e^{-x}}
=\frac{\e^{-x}-1}{\e^{-x}+1}
=1-\frac{2}{\e^{-x}+1}.
\]

3. There are perovskite solar cell concepts consisting also of organic semiconductor materials \cite{Schmidt2021RoadMap}. 
E.g.\ in \cite{ag-kna19} the fullerene $C_{60}$ was used as electron transport material.
Thus, for the description of  charge transport in organic semiconductor materials,  Gauss--Fermi integrals have to be used for the statistical relation. According to \cite[Subsec.\ 2.1]{ag-gli18}, the Gauss--Fermi integrals fulfil similar essential properties as the Blakemore statistics $F_{B,\gamma}$ for $\gamma=1$  used in the discussion here. The techniques, how to derive positive lower and upper bounds below the number of transport states for the organic species can be found in \cite[proofs of Lemma 4.3, Thm. 5.2]{ag-gli18}.
 In the spirit of 
\cite[Section~6]{ag-gli24} solvability and corresponding bounds for the densities could be realized. Thus, by the methods of the present paper, a uniqueness result is to be expected in this setting, too.


\paragraph{Acknowledgement.}
The authors thank Hannes Meinlschmidt and Joachim Rehberg for fruitful discussions on 
integrability properties of solutions to quasilinear 
parabolic equations.

This work was supported by the Deutsche 
Forschungsgemeinschaft (DFG, German Research Foundation) under 
Germany's Excellence Strategy -- The Berlin Mathematics Research 
Center MATH+ (EXC-2046/1, project ID: 390685689).

\small
\bibliography{pero}

\end{document}